\documentclass[11pt]{amsart}

\usepackage{amssymb, amsmath, amscd, amsthm, color, epsfig}

\usepackage[all]{xy}          
\xyoption{dvips}              

\newcommand{\R}{{\mathbb R}}

\newcommand{\delh}[1]{\Delta_{\text{holim}}^{#1}}
\newcommand{\delc}[1]{\Delta_{\text{conf}}^{#1}}

\newcommand{\tdo}{\mathcal{D}^{o}}
\newcommand{\tde}{\mathcal{D}^{e}}
\newcommand{\Td}{\mathcal{D}}
\newcommand{\G}[2]{\Gamma_{#1,#2}}

\newcommand{\K}{{\mathcal{K}}}
\newcommand{\W}{{\mathcal{W}}}

\newcommand{\FM}[2]{F[#1, #2]}

\newcommand{\HO}{{\mathcal{H}}}

\theoremstyle{plain}
\newtheorem{thm}{Theorem}[section]
\newtheorem{prop}[thm]{Proposition}
\newtheorem{lemma}[thm]{Lemma}

\theoremstyle{definition}
\newtheorem{definition}[thm]{Definition}

\theoremstyle{remark}
\newtheorem{rem}{Remark}

\newcommand{\refT}[1]{Theorem~\ref{T:#1}}

\newcommand{\refP}[1]{Proposition~\ref{P:#1}}
\newcommand{\refD}[1]{Definition~\ref{D:#1}}
\newcommand{\refL}[1]{Lemma~\ref{L:#1}}

\begin{document}

\title{Configuration space integrals and Taylor towers for spaces of 
knots}
\author{Ismar Voli\'c}
\address{Department of Mathematics, University of Virginia, 
Charlottesville, VA}
\email{ismar@virginia.edu}
\urladdr{http://www.people.virginia.edu/\~{}iv2n}

\begin{abstract}
We describe Taylor towers for spaces of knots arising from 
Goodwillie-Weiss calculus of the embedding functor and extend the 
configuration space integrals of Bott and Taubes from spaces of knots 
to the stages of the towers.  
We show that certain 
combinations of integrals, indexed by trivalent diagrams, yield 
cohomology classes of the stages of the tower, just as they do for 
ordinary knots.  
\end{abstract}

\maketitle

{\tableofcontents}

\section{Introduction}\label{S:Intro}

In this paper we use configuration space integrals to
establish a concrete connection 
between the study of knots and Goodwillie-Weiss calculus of the embedding functor \cite{We, GW}.  We do this by factoring 
the Bott-Taubes map, well-known 
to knot theorists, through a tower of spaces arising from
this theory.

In more detail, fix a linear inclusion of $\R$ into $\R^m$.  We study
 \emph{long knots}, namely embeddings of $\R$ in $\R^{m}$ which agree with this 
linear map outside of a compact set.  The space of such knots, with 
compact-open topology, is homotopy 
equivalent to the space of 
\emph{based} knots in $S^{m}$.  These can be thought of as maps of $S^{1}$
 ``anchored'' at, say, 
the north pole, or, as we prefer, maps of the interval $I$ to 
$S^{m}$ which are embeddings except at the endpoints.  The endpoints 
are mapped to the north pole with the same derivative.  It is not hard to see that this space of based knots is a deformation retract of the space of based knots in the sphere which are prescribed in a neighborhood of the north pole.  The latter, on the other hand, is clearly homotopy equivalent to the space of long knots.  Let $\K_{m}$ be 
the space of long knots in 
 $\R^{m}$ or $S^{m}$, $m\geq 3$.  To simplify 
 notation, we will often set $\K=\K_{3}$ when we wish to distinguish 
 the case of classical knots from all others.

  At the heart of our results are Bott-Taubes
  configuration space integrals \cite{BT} which are used for producing cohomology classes on 
 $\K_{m}$.  They were originally defined for ordinary knots, i.e. embeddings of $S^1$ in $\R^3$, but the modification to long knots is straightforward \cite{Catt2}.  The idea is to start with a chord diagram with $2n$ 
vertices joined by chords, evaluate a knot on as many 
points, and then consider $n$ 
maps to spheres given by normalized differences of 
pairs of those points.  Which points are paired off is prescribed by 
the chord diagram.  Pulling back the product of volume forms on the 
spheres via a product of these maps yields a form on the product of $\K_{m}$ with a suitably compactified configuration space of $2n$ points in 
$\R^{m}$.  This form can then be pushed forward to $\K_{m}$.  Various arguments involving Stokes' Theorem and the combinatorics of chord diagrams in the end guarantee that the result is a cohomology class.  This was first done by Altschuler and Freidel 
\cite{Alt} and D. Thurston \cite{Th} for 
$\K$, and then generalized by Cattaneo, Cotta-Ramusino, and Longoni 
\cite{Catt} to $\K_{m}$.  We will recall the main features of Bott-Taubes integration in \S\ref{S:B-TIntegrals}.  We will not provide all the details since they can be found in D. Thurston's work \cite{Th} or the survey paper \cite{Vo3}.

The other ingredient we need is the Taylor tower for $\K_{m}$ arising 
from the calculus of the embedding functor.  One considers spaces of ``punctured knots," or embeddings of the interval with some number of subintervals removed.  These spaces fit into cubical diagrams whose homotopy limits define stages of the tower, or ``Taylor approximations" to $\K_{m}$.  For $m>3$, the tower converges (see \refT{Connectivity} for the precise statement) so it represents a good substitute for $\K_m$.  We review the construction of the tower in some detail in \S\ref{S:G-WConstruction}.  Since embedding calculus is the less familiar half of the background we require, we do not assume the reader has had previous exposure to it.  

In \S\ref{S:Graph}, we then turn our attention to extending the Bott-Taubes integrals to the tower and deduce our main result, stated more precisely as \refT{MainTheorem}.

\begin{thm}\label{T:IntroMain}  Bott-Taubes integrals factor through the stages of the Taylor tower for $\K_{m}$, 
$m\geq 3$.
\end{thm}

One importance of this theorem is that the stages of the Taylor tower lend themselves to a geometric analysis which complements the combinatorics and integration techniques of Bott and Taubes.  In particular, one might ask if \emph{all} cohomology classes of spaces of knots arise through Bott-Taubes integration and proceed to look for the answer in the Taylor tower.  Something along these lines has been done for the case of classical knots $\K$ where some, but not all, of the constructions and results presented here hold as well.  In particular, Bott-Taubes integration produces knot invariants and it was shown in \cite{Vo} that the Taylor tower for $\K$ in fact classifies finite type (Vassiliev) invariants.  \refT{IntroMain} plays a crucial in establishing this result.   The hope is that examining the Taylor tower more closely will shed new light on finite type invariants and the slightly mysterious appearance of integration techniques in knot theory. Some more details will be given at the end.

\subsection{Acknowledgements}  I am grateful to Tom 
Goodwillie for the guidance he provided and knowledge he shared with 
me over the years.  I am also endebted
 to Riccardo Longoni for his help 
 with the combinatorics of chord diagrams, as well as to Dev Sinha, Pascal Lambrechts, and especially Greg Arone for comments and suggestions.

\section{Goodwillie-Weiss construction of the Taylor tower for $\K_m$}\label{S:G-WConstruction}

Let $M$ and $N$ be smooth manifolds of dimensions $m$ and $n$, and 
let $Emb(M,N)$
denote the space of embeddings of $M$ in $N$.  Weiss \cite{We} (see 
also \cite{GKW})
develops a certain tower of spaces for studying $Emb(M,N)$.  Its stages $T_{r}$ 
 are constructed from spaces of embeddings of some simple 
codimension $0$ submanifolds of $M$ in $N$.  Each $T_{r}$ comes with a 
canonical map 
from $Emb(M,N)$ and to $T_{r-1}$,
and is in principle 
easier to understand than $Emb(M,N)$ itself.  Goodwillie and Weiss 
\cite{GW}, using work of Goodwillie and Klein \cite{GK}, then prove 
the following

\begin{thm}\label{T:Connectivity}
If 
$n\!-\!m\!>\!2$, the map $Emb(M,N)\longrightarrow T_{r}$ is 
$(r(n-m-2)+1-m)$-connected.
\end{thm}
Since the connectivity increases with $r$,
the inverse limit of the tower is weakly equivalent to $Emb(M,N)$.
Spaces $T_{r}$ are examples of ``polynomial,'' or ``Taylor,'' approximations of 
$Emb(M,N)$ in the sense of Goodwillie calculus.

The general definition of the stages of the Taylor tower can be found in \cite[Section 5]{We}.  However, in the case of $\K_m$, the definition readily simplifies to a concrete construction which produces an equivalent tower, even for classical knots (the edge of the dimensional assumption in the above theorem) \cite[Section 5.1]{GKW}.  We thus focus in some detail on the construction of the Taylor tower for spaces of knots and start with some general definitions.

\begin{definition}
A \emph{subcubical diagram} $C_{r}$ is a functor from 
the category of nonempty subsets $S$ of $\{1, \ldots, r\}$ with inclusions as morhisms to 
spaces,
i.e. it is a diagram of $2^{r}-1$ 
spaces $X_{S}$ so that, for every containment 
$S\subset S\cup \{i\}$, there is a map $X_{S}\to X_{S\cup \{i\}}$ and 
every square 
$$\xymatrix{
X_{S} \ar[r]  \ar[d]  &  X_{S\cup\{i\}} \ar[d] \\
X_{S\cup\{j\}} \ar[r] &  X_{S\cup\{i,j\}}
}
$$
commutes.
\end{definition}
Now let $x_{1}, 
\ldots, x_{r}$, be the barycentric coordinates of the standard $(r-1)$-simplex, which we denote by $\delh{r-1}$.  
Denote by $\delh{S}$ the face of $\delh{r-1}$ given by 
$
x_{i}=0 \ \ \text{for all}\ \ i\notin S.
$
Thus if $T\subset S$, we have an inclusion 
$
\delh{T}\hookrightarrow \delh{S}
$
of a particular face of $\delh{S}$.

\begin{definition}\label{D:CubeHolim}
The \emph{homotopy limit} of an $r$-subcubical diagram 
$C_{r}$, denoted by $holim(C_{r})$, is a subspace of the space of smooth maps
$$
\prod_{\emptyset\neq S\subseteq\{1, \ldots, r\}} Maps(\delh{S}, 
X_{S})
$$
consisting of collections of smooth maps $\{\alpha_{S}\}$ such that, for 
every map $X_{S}\to X_{S\cup \{i\}}$ in the diagram, the square
$$\xymatrix{
\delh{S}\ar[r]^{\alpha_{S}} \ar@{^{(}->}[d] & X_{S} \ar[d] \\
\delh{S\cup\{i\}} \ar[r]^{\alpha_{S\cup \{i\}}}   &   X_{S\cup \{i\}}
}
$$
commutes.
\end{definition}

\begin{rem}
We will want to define certain forms on our homotopy limits in \S\ref{S:Graph} so we consider only smooth maps in the above definition, thereby obtaining differentiable spaces.  If we had instead considered spaces of all maps from simplices, we would have obtained homotopy equivalent spaces.  More on homotopy limits of diagrams in model categories can be found in \cite{BK, DS}.
\end{rem}

Since $C_{r}$ contains $C_{r-1}$, there are  
projections
$
holim(C_{r})\to
holim(C_{r-1})
$
for all $r>1$.
Further, if $X_{\emptyset}$ fits $C_{r}$ as its initial 
space, i.e. it maps to all other spaces in $C_{r}$ 
and makes all the resulting squares commutative (and hence it maps to $holim(C_{r})$), the diagram
\begin{equation}\label{E:CommutativeTriangle}
\xymatrix{
X_{\emptyset} \ar[r] \ar[dr] & holim(C_{r}) \ar[d] \\
                          & holim(C_{r-1})
}
\end{equation}
commutes.

\noindent
We can now define the Taylor tower for the space of 
knots.
For $r>1$, let
$\{A_{i}\}$, $1\!\leq\! i\!\leq r$, be a collection of 
disjoint 
closed subintervals of $I\subset \R$,
indexed cyclically.  
For each nonempty subset $S$ of \{1, \ldots, r\}, define the space of 
maps
$$
E_{S}= \mbox{Emb}(I\!-\! \bigcup_{i\in S}A_{i},S^{m})
$$ 
which are smooth embeddings other than at the endpoints of $I$.  The 
endpoints are, as usual, mapped to the north pole in $S^{m}$ with the 
same derivative.

The $E_{S}$ can be thought of as spaces of ``punctured knots,''  
and are
path-connected even for $m=3$  since any punctured knot can be isotoped to the 
punctured unknot by ``moving strands through the holes''.
If $T\subset S$, there is a restriction $E_{T}\to E_{S}$
which simply sends a punctured knot to the same knot with more 
punctures.  These restrictions clearly commute.  We can thus make the following

\begin{definition}\label{D:SubcubicalDiagram}
Denote by $EC_{r}$ the subcubical diagram sending $S$ to $E_S$ for all nonempty seubsets $S$ of $\{1, ..., r\}$ and sending inclusions to restrictions.
\end{definition}

The homotopy limit of this diagram is the central object of study here so we give some details about what \refD{CubeHolim} means in this case.  Keeping in mind that a path in a space 
of embeddings is an isotopy, a point in $holim(EC_{r})$ is a list of 
embeddings and families of isotopies:
\begin{itemize}
\item an embedding $e_{i}\in E_{\{i\}}$ for each $i$;

\item an isotopy $\alpha_{ij}\colon\delh{1}\to E_{\{i,j\}}$ for each $\{i,j\}$ such that 
$$\alpha_{ij}(0)=e_{i}\vert_{_{E_{\{i,j\}}}}, \ \  
\alpha_{ij}(1)=e_{j}\vert_{_{E_{\{i,j\}}}};$$
\item a 2-parameter isotopy $\alpha_{ijk}\colon \delh{2}\to E_{\{i,j,k\}}$ for each ${\{i,j,k\}}$
whose restrictions to the faces of $\delh{2}$ are
$$\alpha_{ij}\vert_{_{E_{\{i,j,k\}}\times\delh{1}}},\ \ 
\alpha_{jk}\vert_{_{E_{\{i,j,k\}}\times\delh{1}}},\ \ 
\alpha_{ik}\vert_{_{E_{\{i,j,k\}}\times\delh{1}}};
$$
and in general,
\item each $(|S|-1)$-parameter isotopy 
$\delh{|S|-1}\to E_{S}$ is 
determined on the face of $\delh{|S|-1}$ by the restriction of a $(|S|-2)$-parameter 
isotopy of a knot with $|S|-1$ punctures to the same isotopy of a knot with 
one more puncture.
\end{itemize}

Since we chose a definite indexing for the subintervals $A_{i}$ of $I$, $i\in\{1, 
\ldots,r\}$,  
and thus for spaces of punctured knots $E_{S}$, 
$S\subseteq\{1,\ldots,r\}$, there are canonical maps 
$$holim(EC_{r})\longrightarrow holim(EC_{r-1}),\ \ \ r>2.$$

Also, $\K_{m}$ maps to each 
$E_{S}$ again by restriction.   Every square face in the cubical diagram obtained by adjoining $\K_{m}$ in the missing corner of $EC_{r}$ 
commutes, so that we get commutative triangles as in \eqref{E:CommutativeTriangle}.  
 
\begin{definition}\label{D:Tower}  For all 
$r>0$, let $\HO_{r}=holim(EC_{r+1})$ be the \emph{$r$th stage} of the \emph{Taylor tower 
for the space of knots},
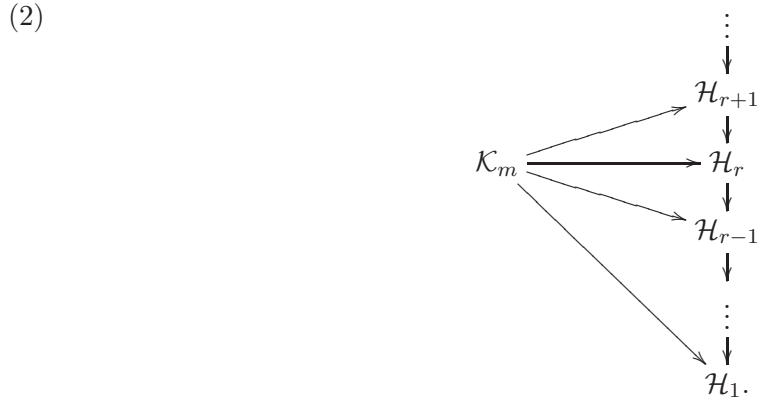
\begin{equation}\label{E:GoodwillieTower}
\xymatrix@R=10pt@C=60pt{
         &    \vdots \ar[d]   \\
         &      \HO_{r+1} \ar[d] \\
   \K_{m} \ar[r] \ar[dr] \ar[dddr]  \ar[ur] &  \HO_{r}\ar[d] \\
        &     \HO_{r-1}  \ar[d]  \\
        &    \vdots \ar[d]  \\
        &     \HO_{1}.
}
\end{equation}
\end{definition}
The tower is shown here with the canonical maps from $\K_{m}$.  Note that \refT{Connectivity} implies convergence of this tower  to $\K_m$ as long as $m>3$.

\begin{rem}
Each point in $EC_r$ determines a knot as 
long as $r>2$.  In fact, we only need to know what the elements of 
such a compatible collection are in $E_{\{1\}}, \ldots, E_{\{r\}}$ in 
order to recover a knot.
Thus $\K_{m}$ actually completes the subcubical diagram of punctured knots as 
its limit for $r>2$. We are therefore in some sense attempting to understand $\K_{m}$, a limit of a certain diagram, 
by instead studying its homotopy limit.
\end{rem}

Spaces $\HO_{r}$ are precisely what $\K_{m}$ will be replaced by in 
the Bott-Taubes construction of the next section.

\section{Bott-Taubes configuration space integrals}\label{S:B-TIntegrals}

\subsection{Trivalent diagrams}\label{S:Trivalent}

Before we turn to configuration space integrals, we give a very brief introduction to a class of diagrams which turns 
out to best keep track of the combinatorics associated to those 
integrals.  More details can be found in \cite{BN, Catt, Long}.

\begin{definition}\label{D:TrivalentDiagrams}  Let \emph{trivalent 
diagram of degree $n$} be a connected graph consisting of an oriented interval, $2n$ vertices, and some number of chords between them.  The vertices lying on the interval are called \emph{interval} and are connected to the rest of the graph by exactly one chord.  The vertices not on the interval are \emph{free} and have exactly three chords emanating from them.
\end{definition}

Depending on whether we are working in $\R^{m}$ for $m$ even or odd, 
our configuration space integrals may change sign due to a permutation of the configuration points, a permutation in the product of certain maps to spheres, or due to a composition of one of those maps with the antipodal map.  These sign changes correspond to a permutation of 
the vertices or chords of a trivalent diagram, or change in the orientation of a chord (see discussion following \refT{Thurston} for more details).  As in \S4.1 of \cite{Catt}, we thus distinguish two classes of diagrams as follows.

\begin{figure}[h]
\begin{center}
\input{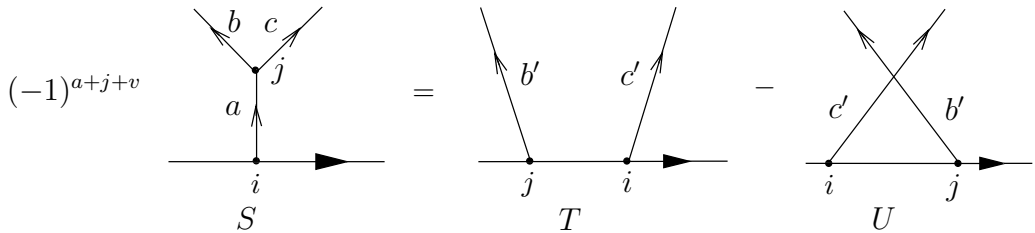}
\caption{$STU_{e}$ relation.  The three diagrams agree outside the 
pictured portions.  Here $v$ is the number of 
interval vertices of the diagram $S$.}\label{F:STU}
\end{center}
\end{figure}

Label the vertices of a trivalent diagram by $1, ..., 2n$, orient its chords, and let $TD_{n}^o$ be the set of all 
trivalent diagrams of degree $n$ with these decorations.  Define $TD_{n}^e$ in the same way except also label the chords.  Let $STU_{e}$ be the relation 
from Figure \ref{F:STU}.  The decorations on the three diagrams in the picture should be compatible:  Since the diagrams are the same outside the pictured portions, the vertex labels and orientations of chords and identical there.  This leaves chord labels.  In the only part where diagrams $S, T$, and $U$ differ, the chords are labeled as in the figure, with $b'=b$ if $b<a$ and $b'=b-1$ if $b>a$.  Same for $c'$.  We follow this pattern outside the pictured parts, and again note that now the chords for $T$ and $U$ are the same as those in $S$.  Thus each chord for $T$ and $U$ is labeled as the corresponding chord in $S$ unless its label is greater than $a$, in which case it is decreased by one.  

Finally let $STU_{o}$ be the same relation as $STU_{e}$
except the factor of $(-1)^{a+j+v}$ is taken away, as are all the chord labels.

\begin{definition} Let $\tdo_{n}$ and $\tde_{n}$ be real vector spaces 
generated by $TD_{n}^o$ and $TD_{n}^e$, modulo the $STU_{o}$ and $STU_{e}$ relations, 
respectively, with
\begin{itemize}
\item Diagrams containing a chord connecting two consecutive interval vertices, diagrams containing a double chord, and diagrams connecting a vertex to itself are all set to zero.
\item For $D_{1}, D_{2}\in \tdo_{n}$ which differ in the 
 orientation of chords, set 
$D_{1}=(-1)^{s}D_{2}$, where $s$ is the number of chords with at least one free end vertex whose orientation 
is different.

\item For $D_{1}, D_{2}\in \tde_{n}$ which differ in the 
 orientation and labels of chords, set 
$D_{1}=(-1)^{s}D_{2}$, where $s$ is sum of the number of chords with at least one free end vertex whose orientation 
is different and the order of the permutation 
of the chords.
\end{itemize}
\end{definition}


\begin{rem}
The relation which sets a diagram containing a chord connecting two consecutive interval vertices to zero is usually called the $1T$ ({\em one-term}) relation, and it is taken away if one considers {\em framed} knots.  It is also intimately related to the correction term $M_DI(D_1,K)$ appearing in \refT{Thurston} and \refT{MainTheorem} \cite{BT, Th}.


\end{rem}

 Let $\tdo=\oplus_{n>0}\tdo_{n}$ and $\tde=\oplus_{n>0}\tde_{n}$.  Since our arguments do not depend of which space of diagrams is considered, we 
will just let $\Td$ stand for either from now on and make  
some remarks on the parity where needed.  Same for $TD_{n}^o$ and $TD_{n}^e$ which we will denote by $TD_{n}$.  

\begin{definition}
Let $\W$ be the space of \emph{weight systems} defined as the dual of 
$\Td$.  Let $\W_{n}$ be the degree $n$ part of $\W$.  
\end{definition}

\begin{thm}[\cite{BN}, Theorem 7] $\Td$ and $\W$ are Hopf algebras.
\end{thm}

The product on $\Td$ is given by continuing the interval of one diagram into another, and the coproduct is essentially given by breaking up the diagram into connected pieces (see Definition 3.7 of \cite{BN}).  A consequence of the theorem is that it suffices to consider only \emph{primitive} weight systems, as we will do from now on.  These are precisely the weight systems which vanish on products of diagrams \cite{BN}.

\subsection{Integrals and cohomology 
classes}\label{S:IntegralsCohomology}

Recall that the linking number of two knots can be obtained by taking two points, one on each knot, and integrating over $S^1\times S^1$ the pullback of the volume form on $S^2$ via the map giving the direction between those two points.  Bott-Taubes configuration space integrals are in a way generalizations of this procedure to the case of a single knot.  However, the points could now collide, so this configuration space has to be compactified for integration to make sense.  Thus given a smooth manifold $M$ of dimension $m$, let $F(k,M)$ be the configuration space of $k$ distinct points in $M$ and let $F[k,M]$ be its \emph{Fulton-MacPherson compactification} \cite{FM, AS}.  The standard way to define this space is through blowups of all the diagonals in $M^k$, but an alternative definition which does not use blowups was given by Sinha.  We state it here in the relevant case of $M=\R^m$. 
\begin{definition}[\cite{Dev1}, Definition 1.3]
Let $F[k,\R^m]$ be the closure of the image of $F(k,\R^m)$ in $(\R^{m})^k\times (S^{m-1})^{k\choose 2}\times [0,\infty]^{k\choose 3}$ under the map which is the inclusion on the first factor and on the second and third sends the point $(x_1, ..., x_k)$ to the product of all $\frac{x_i-x_j}{|x_i-x_j|}$ and $\frac{|x_i-x_j|}{|x_i-x_l|}$, $1\leq i<j<l\leq k$, respectively.
\end{definition}
The compactification $F[k,\R^m]$ is a smooth manifold with corners of dimension $km$ \cite[\S3]{Dev1}, i.e. a space whose every point has a neighborhood homeomorphic to
$$
\R^d\times \R_{+}^{km-d}
$$ 
for some $d$ and such that each transition function extends to an embedding of a neighborhood containing its domain.  It is also compact in the more general case when $M$ is compact.  The configuration points in $F[k,\R^m]$ are allowed to come together while the directions as well as the relative rates of approach of the colliding points are kept track of. Codimension one faces (strata, screens), important for Stokes' Theorem arguments, are given by some number of points colliding at the same time.  The combinatorics of these compactifications are very interesting and deep, and have been related to Stasheff associahedra and certain spaces of trees \cite[\S4]{Dev}.

To make Stokes' Theorem arguments work out, Bott and Taubes make the following definition.

\begin{definition}[\cite{BT}, page 5283]\label{D:BTPullback}
Define $F[k,s; \K_{m}, S^{m}]$ as the pullback of 
$$
\xymatrix{
  \FM{k}{I}\times \K_m  \ar[rr]^-{\text{evaluation}} &   &   \FM{k}{S^m}   & &
 \FM{k+s}{S^m}  \ar[ll]_-{\text{projection}}.
}
$$
\end{definition}

These spaces are suitable for integration, as we have

\begin{prop}[\cite{BT}, Proposition A.3]\label{P:BTBundle}
Spaces $F[k,s; \K_{m}, S^{m}]$ fiber over $\K_m$ and the fibers are
smooth compact manifolds with corners.
\end{prop}

Each fiber of $F[k,s; \K_{m}, S^{m}]$ over $\K_m$ can be thought of as a configuration space of $k+s$ points in $S^m$ with $k$ of them constrained to lie on some knot $K\in\K_m$.  The connection to trivalent diagrams is now clearer; the configuration points which can be anywhere in $S^m$ can be represented by the free vertices while those which lie on a knot can be represented by the interval ones.

Since we wish to consider directions between points, we replace $S^{m}$ by 
$\R^{m}\cup\infty$.  This in turn replaces based knots in $S^{m}$ by long knots 
in $\R^{m}$, but introduces ``faces at infinity" discussed after \refT{Thurston} and in \refL{KnotInfinityFaces}.

Now suppose a 
labeled trivalent diagram $D\in TD_{n}$ with $k$ 
interval and $s$ free vertices is given (so $k+s=2n$).  A chord connecting vertices $i$ and $j$ gives a map 
\begin{align}\label{E:PullbackMap}
h_{ij}\colon F[k,s;\K_m,\R^{m}] & \longrightarrow S^{m-1}  \\
(p_1, ..., p_i, ..., p_j, ..., p_{k+s}) & \longmapsto
\frac{p_j-p_i}{|p_j-p_i|}.\notag
\end{align}
The product of these maps over all $(k+3s)/2$ chords of $D$ can be used for pulling back the product of unit volume forms $\omega_{ij}$, call it $\omega$, from the product of spheres $S^{m-1}$ to $F[k,s;\K_m,\R^{m}]$.
We denote the resulting $(k+3s)(m-1)/2$--form on $F[k,s;\K_m,\R^{m}]$ by $\alpha$. Because of \refP{BTBundle}, it makes sense to push this form forward to $\K_m$, i.e. integrate it along the fiber of the map
$$
\pi\colon F[k,s;\K_m,\R^{m}]\longrightarrow \K_m.
$$  
Finally let $I(D,K)$ stand for the pullback of $\omega$ followed by this pushforward $\pi_*\alpha$:
$$
\xymatrix{
\Omega F[k,s;\K_m,\R^{m}]  \ar[d]^{\pi_*}  &  & \Omega (S^{m-1})^{(k+3s)/2} 
\ar[ll]-_{\prod\limits_{\text{chords }ij}h_{ij}^*}  
\ar[dll]^{I(D,K)} \\
\Omega \K_m
}
$$ 
Since the fiber of $\pi$ has dimension $k+ms$, the resulting form on $\K_m$ has dimension
$$
\frac{k+3s}{2}(m-1)-(k+ms)=(m-3)\frac{k+s}{2}=(m-3)n.
$$
This is not necessarily a closed form.  However, let $D_1$ be the diagram consisting of two interval vertices and one chord between them.
We then have
\begin{thm}\label{T:Thurston} 
 For a nontrivial primitive weight system $W\in\W_{n}$, the map
 $
 T(W)\colon\K\to \R
 $
 given by 
$$
K\longmapsto \frac{1}{(2n)!}
\sum\limits_{D\in TD_n}
W(D)(I(D,K)-M_DI(D_1,K)),
$$
represents a nontrivial element of $H^{(m-3)n}(\K_{m};\R)$.
Here $M_D$ is a real number which depends only on $D$ and $M_DI(D_1,K)$ vanishes for $m>3$.
\end{thm}
In the case $m=3$, this theorem was first proved for ordinary closed knots by Thurston \cite{Th} and Altschuler 
and Friedel \cite{Alt} who also show the zeroth cohomology class one gets this way on $\K$ is in fact a finite type $n$ invariant.  The generalization to $m>3$
is due to Cattaneo, Cotta-Ramusino, and Longoni, who also show the cohomology classes obtained this way are nontrivial \cite[Section 6]{Catt}. 
The proof does not depend on $m$ except a little care has to be taken with signs.
Since a labeling 
 of a diagram determines the labeling of configuration points in 
 $F[k,s;\K_m,\R^{m}]$, changing the orientation of $D$ may affect the signs of $I(D,K)$ and $M_DI(D_1,K)$ depending on $m$ (orientation of the fiber might change).  But the two diagram algebras $\mathcal{D}^{e}_{n}$ and $\mathcal{D}^{o}_{n}$, corresponding to $m$ even and odd, are defined precisely so that $W$ depends on the sign in the same way.

One proof of \refT{Thurston} is via Stokes' Theorem and proceeds by checking that the integrals on the boundary of the fiber of $\pi$ 
either vanish or cancel out within the sum, so that the sum is in fact a closed form.  Different arguments are 
used for various types of faces, which are called
\emph{principal} if exactly two points 
degenerate; \emph{hidden} if more than two, but not all, points degenerate; and \emph{faces at infinity} if one or more points approach 
infinity.  The correction term $M_DI(D_1,K)$ comes from the possible contribution of the \emph{anomalous} face corresponding to all configuration points coming together \cite[Proposition 4.8]{Vo3}.  While it is easy to see that the contribution is zero in case of knots in $\R^m$, $m>3$ \cite[Proposition 6.3]{Vo3}, it is a conjecture that this is also the case for $m=3$.  D. Thurston \cite{Th} and Poirier \cite{Poir} have computed it to be zero in some simple cases.

The vanishing arguments can be found in \cite{BT, Th, Poir, Vo3} and can be written down very concretely using explicit coordinates on compactified configuration spaces \cite[page 5286]{BT} (see also \S4.1 in \cite{Vo3}).  Integrals along principal faces do not necessarily vanish, but they can be grouped according to the $STU$ relation (and another relation which follows from it, usually called the $IHX$ relation; see \cite{BN, Vo3}).  These integrals then cancel in the sum of \refT{Thurston} \cite[\S4.4]{Vo3}.  For other faces, a key observation time and again is that the product of the maps $h_{ij}$ factors through a space of lower dimension than the product of the spheres which is its initial target.  Therefore $\alpha$ must be zero.  This type of argument is illustrated in \refL{KnotInfinityFaces} below and it immediately takes care of the vanishing of integrals along hidden faces \cite[Proposition 4.4]{Vo3} and faces at infinity  where one or more of the points off the knot go to infinity \cite[Proposition 4.7]{Vo3}.  In case of long knots, however, there is an extra case of such a face corresponding to some points on the knot going to infinity.  This cannot happen in the Bott-Taubes/Thurston setup since they consider closed knots.  We deal with this case in the following lemma.

\begin{lemma}\label{L:KnotInfinityFaces}
The pushforward $I(D,K)$ vanishes on the faces of the fiber of $\pi$ corresponding to some or all points on the knot going to infinity.
\end{lemma}

\begin{proof}
The argument is essentially that of Proposition 4.7 in \cite{Vo3}.  Recall that our long knots are ``flat" outside a compact set, i.e. they agree with a fixed linear inclusion of $\R$ in $\R^m$.  Suppose a point $p_i$ on the knot tend to infinity.  If $p_i$ is related to another point $p_j$ by a map $h_{ij}$ (meaning there is a chord connecting vertices $i$ and $j$ in $D$), then there are four cases to consider.
\begin{enumerate}
\item 
If $p_j$ does not go to infinity, then $h_{ij}$ restricts to a constant map along this face.  The product of all such maps to $(S^{m-1})^{(k+3s)/2}$ then factors through $(S^{m-1})^{((k+3s)/2)-1}$.  The pullback of $\omega$ to $F[k,s;\K_m,\R^{m}]$ thus has to be zero as does $I(D,K)$.
\item If $p_j$ is on the knot and it also goes to infinity (regardless of whether it does so in the same direction as $x_i$), $h_{ij}$ is constant on this face.
\item
If $p_j$ is off the knot and it also goes to infinity, but in a different direction than that of the fixed linear inclusion, $h_{ij}$ is again constant.
\item
If $p_j$ is off the knot and it goes to infinity in the same direction and at the same rate as $p_i$, then $p_j$ is either connected to a point $p_k$ which does not, in which case $h_{jk}$ restricts to a constant map on this face, or it does, in which case we look at all other points $p_k$ is related to by maps.  Since $D$ is connected, there must eventually be two points for which the map restricts to a constant map (if not, this means the entire configuration is translated along the knot to infinity and this is not a face). 
\end{enumerate}
\end{proof}
%


We next modify the construction outlined in this section to the setting of the Taylor tower and generalize \refT{Thurston}.

\section{Generalization to the stages of the Taylor 
tower}\label{S:Graph}

Remember from \S\ref{S:G-WConstruction} that a point $h$ in $\HO_{k}$ 
is a collection of families of embeddings parametrized by simplices of 
dimensions $0, \ldots, k$.  The families are compatible in the sense 
that a $k$-simplex $\delh{k}$ parametrizes a family 
of knots with $k+1$ 
punctures, while each of its faces parametrizes a family of knots 
with fewer punctures (how many and which punctures depends on  
which barycentric coordinates of $\delh{k}$ are 0).  
However, the evaluation of a punctured knot on a point in $\FM{k}{I}$ may 
not be defined since the configuration points may land in the parts 
of $I$ that have been removed.
To get around this, we will devote most of this section to the 
construction of a smooth map  
$$
\FM{k}{I}\to \delh{k}
$$
whose graph will serve the purpose of choosing a punctured 
knot in the family $h\in\HO_{k}$ depending on where the $k$ points in $I$ may be.
\vskip 4pt
\noindent
The interior of $\FM{k}{I}$, the open configuration 
space $F(k,I)$, is given by points $(x_{1}, \ldots, x_{k})$ which satisfy
$
0<x_{1}<x_{2}<\cdots<x_{k}<1.
$
Thus we have a natural identification
$
F(k,I)\simeq \Delta^{k},
$
where $\Delta^{k}$ denotes the open $k$-simplex.
Let $\delc{k}$ be the closed simplex identified 
with the obvious compactification of $F(k,I)$, i.e.\!\! adding the faces 
to $\Delta^{k}$.  
Also let $\partial_{i}\delh{k}$ stand for 
the $i$th face of $\delh{k}$ ($i$th barycentric coordinate is 0), and 
let $A(x)$ index the set of holes in which the configuration $x$ may be.  
In other words,
$$
A(x)=\{i :  x_{j} \in A_{i} \mbox{\ for some \ } j\}.
$$

\begin{prop}\label{P:MapConstruction}
There is a smooth map
$
\gamma^{k}\colon \delc{k} \longrightarrow \delh{k},
$
defined inductively, which depends on the choice of the punctures $A_1, ..., A_r$ in $I$.  Moreover, if 
$
\gamma^{i}\colon \delc{i}\to \delh{k}
$
has been defined for all $i<n$,  then, for $1\leq j\leq n-1$,
$
\gamma^{n}\colon \delc{n}\to \delh{k}
$
satisfies
\newcounter{Lcount}
\begin{list}{\roman{Lcount})}{\usecounter{Lcount}}

\item  \begin{gather*}\label{E:GraphCondition1}
\gamma^{n}(x_{1}, \ldots, x_{j-1}, x_{j}, x_{j}, x_{j+2}, \ldots, 
x_{n})=
\gamma^{n-1}(x_{1}, \ldots, x_{j-1}, x_{j}, x_{j+1}, \ldots, 
x_{n-1}), \\ \notag
\gamma^{n}(0, x_{2}, \ldots, x_{n})=
\gamma^{n-1}(x_{2}, \ldots, x_{n}) \\  \notag
\gamma^{n}(x_{1}, \ldots, x_{n-1}, 1)=
\gamma^{n-1}(x_{1}, \ldots, x_{n-1}); \notag
\end{gather*}

\item  There exists an open 
neighborhood $V$ of $x$ and 
\begin{equation*}\label{E:GraphCondition2}
\gamma^{n}(x')\in \bigcap_{i\in A(x)} \partial_{i}\delh{k} 
\mbox{ \ \ for all } x'\in V.
\end{equation*}

\end{list}

\end{prop}

Conditions i) and ii) are required because of the following:  Let $x=(x_{1}, \ldots,x_{k}),\  0\leq 
x_{1}\leq \cdots \leq x_{k}\leq 1$ parametrize $\delc{k}$. 
The image in $\delc{k}$ of two points coming together in 
$F[k,s;\K_{m},\R^{m}]$ is
$x_{j}=x_{j+1}$, 
$1\leq j\leq k-1$.  This situation translates into the 
pushforward of a certain form along a principal face and we wish for 
integrals like this to cancel due to the $STU$ and $IHX$ relations after 
considering sums over all trivalent diagrams.  
The cancellation will 
only be possible if the integrals corresponding to each 
triple of diagrams have the same value when two points collide.  
However, one of 
the diagrams in the $STU$ relation has fewer interval vertices, i.e. it 
is associated with the space $\FM{k-1}{I}$.
A way to ensure the appropriate integrals over $\FM{k}{I}$ and 
$\FM{k-1}{I}$ are equal is to 
define $\gamma^{k}$ inductively based on the number of points in a 
configuration (keeping in mind that $\delc{0}$ is a point) and to further impose condition i).  The last two equations in i) are required for the integrals along the faces 
given by 
points colliding with the basepoint in $S^{m}$ to cancel out.

As for condition ii), given $t$ 
in $\delh{k}$ and $h$ in $\HO_{k}$, one gets a point $h_{t}$ in $\delh{k}\times\mathcal{H}_{k}$ which is an embedding of the interval with up to $k$ punctures.  As mentioned at the beginning of this section, we 
want the evaluation of $h_{t}$ on a configuration
to be defined for all points $(q, t, h)\in \Gamma_{k}\times 
\mathcal{H}_{k}$.  We therefore need that, whenever
 $x\in \delc{n}$, $t=\gamma^{n}(x)$ is a point in $\delh{k}$ such 
that the corresponding embedding $h_{t}$ is defined for $x$.  So depending on 
where $x$ is in $I$, $\gamma^{n}$ 
will map it to the interior or a face (or intersection of faces) of 
$\delh{k}$ according to 
whether some of the $x_{j}$ are in any of the removed subarcs 
$A_{i}$ for
$1\leq i\leq k$.  Condition ii) ensures this and more as it requires $\gamma^{n}$ to 
map a \emph{neighborhood} of every point $x$ to the same face as $x$ 
itself.  This is needed for the resulting graph to be a smooth manifold with corners.  Note that the intersection in condition ii) is nonempty since there 
is always at least one more hole in the interval than the number of points in a 
configuration (the number of configuration points is $n\leq k$, while the 
number of holes in $\HO_{k}$ is $k+1$).

\begin{proof}[Proof of \refP{MapConstruction}]
Assume smooth maps $\gamma^{0}, \ldots, \gamma^{n-1}$ have been defined
 on faces of $\delc{n}$ and satisfy conditions 
i) and ii) (smoothness 
is needed for Stokes' Theorem).  Then we can extend 
locally to a function $\gamma^{n}$ 
on all of $\delc{n}$.  However, we need to check that there are 
neighborhoods $U_{x}$ for every point $x\in \delc{n}$ so that the 
local extensions $\gamma^{n}_{x}$ match on intersections.

Thus, given $x=(x_{1}, \ldots, x_{n}) \in \delc{n}$ with all $x_{j}$ distinct, pick a neighborhood 
$U_{x}$ of 
$x$ such that $A(x')\subset A(x)$ (no $x'$ in $U_{x}$ gets into holes $x$ 
did not get into).  This can be done since the $A_{i}$ are closed 
subintervals.  Then two 
intersecting neighborhoods are both mapped to the same face in 
$\delh{k}$, and so condition ii) is satisfied on 
intersections in 
this case.  
(Condition i) is vacuous here since we are in the interior of 
$\delc{n}$.)

If $x_{j}=x_{j+1}$, so that $x$ is on a face $\delc{n-1}$ of 
$\delc{n}$, choose a $U_{x}$ so that its boundary in 
 $\delc{n-1}$ is 
contained in the neighborhood $V$ from condition 
ii) for the point 
$x=(x_{1}, \ldots, x_{j}, \ldots, x_{n-1})\in \delc{n-1}$ and the map 
$\gamma^{n-1}_{x}$.  Now $\gamma^{n}_{x}$, extended from $V$, maps the whole 
half-ball $U_{x}$ to the same face in $\delh{k}$, and these match to 
define a function on intersections.  

The preceeding easily generalizes to 
those $x$ on lower-dimensional faces of $\delc{n}$.  If there is more than one 
$j$ for which $x_{j}=x_{j+1}$, choose $U_{x}$ such that, {\it for each $j$}, 
the part of the boundary of $U_{x}$ given by 
$$U_{x}\cap \{x': x'_{j}=x'
_{j+1}\}
$$ 
equals $V$, where $V$ has been by induction determined by $x=(x_{1}, \ldots, x_{j}, 
\ldots, x_{k-1})\in \delc{n-1}$ and the map 
$\gamma^{n-1}_{x}$.

Thus $\gamma^{n}$ can be defined locally.  To define it as a smooth 
function on the whole $n$-simplex, let $\{U_{\alpha}\}$ be a finite 
open cover of $\delc{n}$ 
given by neighborhoods $U_{x}$.  Similarly, 
$$\gamma^{n}_{\alpha}\colon 
U_{\alpha} \to \delh{n}
$$ 
are given by the maps 
$\gamma^{n}_{x}$.  Let 
$$
\mu_{\alpha}: U_{\alpha} \to I, \ supp(\mu_{\alpha}) \subset 
U_{\alpha}, \ \sum_{\alpha} \mu_{\alpha}=1, \ \mu_{\alpha}>0,$$
be a partition of unity subordinate to the cover $\{U_{\alpha}\}$,
and note that if two functions $\gamma^{n}_{\alpha}$ and 
$\gamma^{n}_{\beta}$ both satisfy conditions i) 
and ii) on 
$U_{\alpha} \cap U_{\beta}$, so will their average, where averaging is 
done by the partition of unity. 
Thus setting $$\gamma^{m}=\sum_{\alpha}\mu_{\alpha}\gamma^{n}_{\alpha}$$
produces a smooth map from the closed simplex $\delc{n}$ to $\delh{k}$ satisfying 
i) and ii). 
\end{proof}

\begin{rem}\label{R:RemarkOnConstruction}
Instead of using $\HO_{k}$ in this construction,
$\HO_{j}$ could have been used, for any $j> k$.  Then $\delh{j}$ would parametrize a 
family of embeddings in $\HO_{j}$, but we would only be interested in the 
subfamily parametrized by the face $\delh{k}$.  There is no ambiguity 
as to which face is meant since the maps $\HO_{j}\to\HO_{k}$ are 
well-defined.
\end{rem}

The space over which generalized Bott-Taubes integration will take place is now easy to define.  Noting that there is a map 
$$
f\colon \FM{k}{I} \longrightarrow \delc{k}
$$
 which is identity on the 
interior of $\FM{k}{I}$ and forgets the extra information about the 
relative rates of approach of the colliding points, we have

\begin{definition}
Let 
$$\Gamma_{k}=\{(x,t) : t=\gamma^{k}(f(x))\}\subset \FM{k}{I}\times \delh{k}$$
be the graph of the composition
$$
\FM{k}{I} \stackrel{f}{\longrightarrow} \delc{k} \stackrel{\gamma^{k}}
{\longrightarrow} \delh{k}.
$$
\end{definition}  
 
Since $\FM{k}{I}$ and $\delh{k}$ are manifolds with corners, it follows from our construction of $\gamma^k$ that
$\Gamma_{n}$ is a manifold with corners for 
all $n\leq k$.
The generalization of the Bott-Taubes setup from the previous section is now 
straighforward.  In analogy with \refD{BTPullback}, we have

\begin{definition}\label{D:Pullback}
Define $\G{k}{s}$ as the 
pullback
$$
\xymatrix{\Gamma_{k,s} \ar[r] \ar[d] & F[k+s, \R^{m}]\ar[d] \\ 
\Gamma_{k}\times \HO_{k+s} \ar[r] & F[k, \R^{m}]. }
$$
\end{definition}
\begin{rem}\label{R:RestrictionRemark}
Recall that for a point in the homotopy limit coming from a knot, all isotopies are constant.  The manifold $\G{k}{s}$ in this case is therefore precisely $F[k,s; \K_{m}, \R^{m}]$ from the Bott-Taubes setup (and an even more special case is $\G{k}{0}=\Gamma_{k}=F[k,I]$).
\end{rem}

Bott and Taubes' proof of \refP{BTBundle}, which they carry out in a very general setting, applies in our case, so that we immediately get an analogous statement

\begin{prop}\label{P:AnalogyToBTBundle}
Spaces $\Gamma_{k,s}$ fiber over $\HO_{k+s}$ and the fibers are
smooth manifolds with corners.
\end{prop}

With Remark \ref{R:RemarkOnConstruction} in mind, we have chosen to construct $\G{k}{s}$ as 
a bundle over $\HO_{k+s}$ (see comment immediately following 
\refT{MainTheorem} for the reason why).  We have also replaced $S^{m}$ by $\R^{m}$ as before.

The fiber of the map $\G{k}{s}\to \HO_{k+s}$ can now thought of as follows:  Recall that a point in $\HO_{k+s}$ is parametrized by $\delh{k+s}$.  Given $h\in \HO_{k+s}$ and depending on where the points of $F[k,I]$ are, a certain point $t\in \delh{k+s}$ is chosen according to our construction.  This gives a particular punctured knot $h_t$.  A point in the fiber is then a configuration space of $k+s$ points in $\R^m$ with $k$ of them constrained to lie on some punctured knot $h_t$.  Note that the only difference between this and Bott-Taubes setup is a genuine knot $K\in \K_m$ is replaced by the punctured knots $h_t$.

Again given a trivalent diagram $D\in TD_{n}$ with $k$ interval and $s$ 
free vertices, there is a map
\begin{equation}\label{E:hD}
\Big(\prod_{\text{chords }ij}h_{ij}\Big)\colon
\G{k}{s}\longrightarrow (S^{m-1})^{(k+3s)/2}
\end{equation}
given by normalized differences of those pairs of points in $\G{k}{s}$ which 
correspond to pairs of vertices connected by chords in $D$.
Each $h_{ij}$ pulls back the volume form $\omega_{ij}$ to an $(m-1)$--form
$\alpha_{ij}$ on $\Gamma_{k,s}$.  The product of the $\alpha_{ij}$ is 
then a
$(k+3s)(m-1)/2$--form $\alpha$ which can be pushed forward along the 
$(k+ms)$--dimensional fiber to produce an $(m-3)n$--form on $\HO_{k+s}$.  This time we denote the pullback followed by pushforward by $I(D,h)$.  Let $M_DI(D_1,h)$ again be the correction term associated with the 
collision of all points in $\G{k}{s}$, 
so that we may state
a generalization of \refT{Thurston}:
\begin{thm}\label{T:MainTheorem} For a nontrivial primitive weight system $W\in\W_{n}$, 
$n>1$, the map
$
T(W)\colon \HO_{2n} \to \R
$
given by
\begin{equation}\label{E:MainMap}
h \longmapsto \frac{1}{(2n)!}
\sum_{D\in TD_n}W(D)(I(D,h)-M_DI(D_1,h))
\end{equation}
represents a nontrival element of $H^{(m-3)n}(\HO_{2n};\R)$.  The real number $M_D$ again depends only on $D$ and the correction term $M_DI(D_1,h)$
is zero for $m>3$. If 
$h\in\HO_{2n}$ comes from a knot, this is the usual Bott-Taubes map from \refT{Thurston}. 
\end{thm}
\noindent
Note that this is a restatement of \refT{IntroMain}.  Also note that 
the degree $(m-3)n$ is in the range given by \refT{Connectivity} for $m>3$.  
It should now also be clear why $\G{k}{s}$ was defined as a bundle over $\HO_{k+s}$.  Each trivalent diagram in the sum has a total 
of $k+s=2n$ vertices.  The domain of $T(W)$ should be the same 
space regardless of what $k$ and $s$ are.  The proper space to 
define $T(W)$ on is thus $\HO_{2n}$, since $D$ can in the extreme 
case be a chord 
diagram with $s=0$ and $k=2n$.

The main point in \refT{MainTheorem} is that the Bott-Taubes map factors through the Taylor tower.  Spaces $\Gamma_{k,s}$ have been constructed so that this is immediate (see Remark \ref{R:RestrictionRemark}).  To prove that the form on $\HO_{2n}$ given by the map \eqref{E:MainMap} is closed, one can repeat verbatim the arguments given in \cite{BT, Th, Vo3} proving that the form from \refT{Thurston} is closed.  Since these arguments are lengthy but straightforward, we will not repeat them here.  It suffices to say that the main reason why the arguments stay the same is that one can write down coordinates on $\Gamma_{k,s}$ in exactly the same way Bott and Taubes do on $F[k,s;\K_m, \R^m]$.  
These coordinates are for example given in equations (12) of of \cite{Vo3}.  Everything in those equations stays the same except a knot $K$ is replaced by a punctured knot $h_t$, as was already hinted at in the discussion following \refP{AnalogyToBTBundle}.  But all the Stokes' Theorem arguments are based on these coordinates so that \S4.2--4.6 in \cite{Vo3}, where the vanishing results are proved, can now be repeated in exactly the same way.  Since everything therefore immediately carries over from the setting of \refT{Thurston} to ours, it follows that the form given by \eqref{E:MainMap} is closed.

\vskip 6pt
\noindent
To conclude, we briefly indicate how the extension 
of Bott-Taubes integration to the Taylor tower gives another 
point of view on finite type knot theory \cite{BN, BN2}.  The fact that 
configuration space integrals can be used to construct the universal 
finite type knot invariant has been known for some time \cite{Alt, 
Th}.

Note that Bott-Taubes integrals produce 0-dimensional cohomology classes, or knot invariants, in case of $\K$.  Also recall that the Taylor stages $\HO_{k}$ for $\K$ is defined the same way as for $\K_m$, $m>3$.  Let $\HO_{k}^*$ be an algebraic analog for the Taylor stage, obtained by replacing the spaces of punctured knots by cochains on those spaces and taking the algebraic homotopy colimit of the resulting subcubical diagram (this colimit is the total complex of a certain double complex).  One then has canonical maps
$$
H^0(\K)\longleftarrow H^0(\HO_{2n})\longleftarrow H^0(\HO_{2n}^*),
$$
neither of which is necessarily an equivalence (the first because one no longer has \refT{Connectivity}). However,
we have
\begin{thm}[\cite{Vo}, Theorem 6.10]\label{T:MainThesisTheorem}
$H^0(\HO_{2n}^*)$ is isomorphic to the set of finite type $n$ knot invariants.
\end{thm}
It is also not hard to see that one has isomorphisms between $H^0(\HO_{2n}^*)$ and $H^0(\HO_{2n+1}^*)$ \cite[equation (34)]{Vo} so that all the stages of the algebraic Taylor tower are accounted for.  Thus its invariants are precisely the finite 
type invariants.  Configuration space integrals and \refT{MainTheorem} are 
central to the proof of \refT{MainThesisTheorem} but the isomorphism itself is 
given by a simple map based on evaluation of a knot on 
some points.

The Taylor tower is thus a potentially a rich source of 
information about finite type theory.  One interesting question is whether the usual Taylor stages $\HO_{2n}$ contain more than just the finite type invariants.  This issue is very closely related to the conjecture that finite type invariants separate knots.
Some further questions are posed in \S6.5 of \cite{Vo}.


\begin{thebibliography}{99}

\bibitem{Alt} D. Altschuler and L. Freidel.
\ On universal Vassiliev invariants.
\ \emph{Comm. Math. Phys.,}
\ 170(1):~41--62, 1995.

\bibitem{AS} S. Axelrod and I. Singer.   
\ Chern-Simons perturbation theory, II.
\ \emph{J. Differential Geom.,}
\ 39(1):~173--213, 1994.

\bibitem{BN} D. Bar-Natan.  
\ On the Vassiliev knot invariants. 
\ \emph{Topology,} 
\ 34:~423--472, 1995.

\bibitem{BN3} D. Bar-Natan.
\ \emph{Perturbative Aspects of the Chern-Simons Topological Quantum 
        Field Theory.}
\ Ph.D. Thesis, Princeton University, 1991.

\bibitem{BN2} D. Bar-Natan and A. Stoimenow.  
\ The fundamental theorem of Vassiliev invariants.  
\ In Lecture Notes in Pure and Appl. Math., 
\ Vol. 184, 1997.


\bibitem{BT} R. Bott and C. Taubes.  
\ On the self-linking of knots.  
\ \emph{J. Math. Phys.,} 
\ 35(10):~5247--5287, 1994.


\bibitem{BK} A. Bousfield and D. Kan.
\ \emph{Homotopy Limits, Completions, and Localizations.}
\ Lecture Notes in Mathematics, Vol. 304, 1972.


\bibitem{Catt2} A. Cattaneo, P. Cotta-Ramusino, and R. Longoni.  
\ Algebraic structures on graph cohomology.  
\ To appear in {\em J. Knot Theor. Ramif.}.

\bibitem{Catt} A. Cattaneo, P. Cotta-Ramusino, and R. Longoni.  
\ Configuration spaces and Vassiliev classes in any dimension.  
\ \emph{Algeb. Geom. Topol.,} 
\ 2:~949--1000, 2002.
 



\bibitem{DS} W. Dwyer and J. Spalinski.
\ Homotopy theories and model categories.
\ In \emph{Handbook of Algebraic Topology},
\ (I. James, ed.), Elsevier, 1995.

\bibitem{FM} W. Fulton and R. MacPherson.  
\ Compactification of configuration spaces.  
\ \emph{Ann. of Math.,}
\ 139:~183--225, 1994.




\bibitem{GK} T. Goodwillie and J. Klein.  
\ Excision statements for spaces of embeddings.  
\ In preparation.

\bibitem{GKW} T. Goodwillie, J. Klein, and M. Weiss.  
\ Spaces of smooth embeddings, disjunction and surgery.  
\ \emph{Surveys on Surgery Theory,} 
\ 2:~221-284, 2001.

\bibitem{GW} T. Goodwillie and M. Weiss. 
\ Embeddings from the point of view of immersion theory II.
\ \emph{Geom. Topol.,}
\ 3:~103--118, 1999. 












\bibitem{Long} R. Longoni.
\ Quantum field theories and graph cohomology.
\ In preparation.


\bibitem{Poir} S. Poirier.
\ The configuration space integrals for links and tangles in $\R^{3}$.
\ \emph{Algebr. Geom. Topol.}
\ 2:~1001-1050, 2002.

\bibitem{Dev1} D. Sinha.
\ Manifold-theoretic compactifications of configuration spaces.
\ To appear in \emph{Selecta Math.}

\bibitem{Dev} D. Sinha.
\ The topology of spaces of knots.
\ Submitted.  
\ math.AT/0202287, 2005.

\bibitem{Th} D. Thurston. 
\ \emph{Integral Expressions for the Vassiliev Knot Invariants.} 
\ Senior Thesis, Harvard University, 1995.
\ math.QA/9901110.


\bibitem{V} V. Vassiliev.
\ Cohomology of knot spaces.
\ In \emph{Theory of Singularities and Its Applications,}
\ (V. I. Arnold, ed.), Amer. Math. Soc., Providence, 1992.

\bibitem{Vas}  V. Vassiliev.  {\em Complements of discriminants of 
smooth maps:  topology and applications.} Translations of 
Mathematical Monographs, Vol. 98, American Mathematical Society, 
Providence, RI, 1992.

\bibitem{Vo3} I. Voli\'{c}.
\ A survey of Bott-Taubes integration.
\ To appear in {\em J. Knot Theor. Ramif.}

\bibitem{Vo} I. Voli\'{c}.
\ Finite type knot invariants and calculus of functors.
\ To appear in {\em Compos. Math.}


\bibitem{We} M. Weiss. 
\ Embeddings from the point of view of immersion theory I.  
\ \emph{Geom. Topol.,}
\ 3:~67--101, 1999. 

\end{thebibliography}
\end{document}